\documentclass{amsart}
\usepackage{amsmath,amssymb,amsthm,amscd}
\newtheorem{formula}{}[section]
\newtheorem{proposition}[formula]{Proposition}
\newtheorem{corollary}[formula]{Corollary}
\newtheorem{lemma}[formula]{Lemma}
\newtheorem{theorem}[formula]{Theorem}
\theoremstyle{definition}
\newtheorem{definition}[formula]{Definition}
\newtheorem{example}[formula]{Example}
\theoremstyle{remark}
\newtheorem*{remark}{Remark}

\begin{document}

\title
[Cohomology of
solvmanifolds and Morse-Novikov theory]
{Cohomology with local coefficients of 
solvmanifolds and Morse-Novikov theory}
\author{Dmitri V. Millionschikov}
\thanks{Partially supported by                       
the Russian Foundation for Fundamental Research, grant no. 99-01-00090
and PAI-RUSSIE, dossier no. 04495UL}
\subjclass{58A12, 17B30, 17B56 (Primary) 57T15 (Secondary)}
\keywords{Solvmanifolds, nilmanifolds, cohomology, local system,
Morse-Novikov theory, solvable Lie algebras }
\address{Department of Mathematics and Mechanics, Moscow
State University, 119899 Moscow, RUSSIA}
\curraddr{Universit\'e Louis Pasteur, UFR de Math\'ematique et d'Informatique, 7 rue Ren\'e Descartes - 67084 Strasbourg Cedex (France)}
\email{million@mech.math.msu.su}

\begin{abstract}
We study the cohomology
$H^*_{\lambda \omega}(G/\Gamma, {\mathbb C})$ of the deRham complex
$\Lambda^*(G/\Gamma)\otimes{\mathbb C}$ of a compact solvmanifold $G/\Gamma$
with a deformed differential $d_{\lambda \omega}=d + \lambda\omega$, where
$\omega$ is a closed $1$-form.
This cohomology naturally arises in the Morse-Novikov theory.
We show that for a solvable Lie group $G$ with a completely solvable
Lie algebra $\mathfrak{g}$ and a cocompact lattice $\Gamma \subset G$
the cohomology
$H^*_{\lambda \omega}(G/\Gamma, {\mathbb C})$ coincides with
the cohomology  $H^*_{\lambda \omega}(\mathfrak{g})$ 
of the Lie algebra $\mathfrak{g}$
associated with the
one-dimensional representation 
$\rho_{\lambda \omega}: \mathfrak{g} \to {\mathbb K},
\rho_{\lambda \omega}(\xi) = \lambda \omega(\xi)$.
Moreover $H^*_{\lambda \omega}(G/\Gamma, {\mathbb C})$
is non-trivial 
if and only if $-\lambda [\omega]$ belongs to the finite subset
$\{0\} \cup \tilde \Omega_{\mathfrak{g}}$ in  $H^1(G/\Gamma, {\mathbb C})$ well defined in terms of $\mathfrak{g}$.
\end{abstract}

\date{}

\maketitle

\section*{Introduction}
In the begining of the 80-th S.P. Novikov  constructed
(~\cite{N1}, ~\cite{N2})
an analogue of the Morse
theory for smooth
closed $1$-forms 
on a compact smooth manifold $M$.
In particular he introduced the Morse-type inequalities
(Novikov's inequalities)
for numbers $m_p(\omega)$ of zeros 
of index $p$ of a closed $1$-form $\omega$ on $M$.
A lot of papers was devoted to this problem in the following
years (see ~\cite{N4} for references). 
In ~\cite{N3}, ~\cite{Pa}  a method 
of obtaining the torsion-free Novikov
inequalities in terms of the deRham complex of manifold was proposed.
This method was based on Witten's approach ~\cite{W} to the Morse theory.
A. Pazhitnov obtained some important results in this direction in ~ \cite{Pa}. 
The cohomology of the deRham complex
$\Lambda^*(M)$ with the deformed differential $d+\lambda \omega$
coincides with the cohomology $H^*_{\rho_{\lambda \omega}}(M,{\mathbb C})$ 
with coefficients in the local system
$\rho_{\lambda \omega}$ of groups ${\mathbb C}$,
$\rho_{\lambda \omega}(\gamma)=\exp{\int_{\gamma} \lambda \omega}$
and for sufficiently large real numbers $\lambda$
we have the following estimate (see ~\cite{Pa}):
$$ m_p(\omega) \ge \dim H^p_{\rho_{\lambda \omega}}(M,{\mathbb C}), \forall p.$$

L. Alania in ~\cite{Al} studied 
$H^*_{\rho_{\lambda \omega}}(M_n,{\mathbb C})$ of a
class of nilmanifolds $M_n$. He proved that
$H^*_{\rho_{\lambda \omega}}(M_n,{\mathbb C})$ is trivial if
$\lambda \omega \ne 0$. The partial answer for the case $\lambda \omega=0$
was obtained by the present author in ~\cite{Mill}.
In both cases the proof was based on the Nomizu theorem ~\cite{Nz} and 
the computations were made in terms of the
corresponding nilpotent Lie algebra $\mathcal{V}_n$.
The starting
point of this article was the intention
to improve the results of ~\cite{Al} in more general situation
considering solvmanifolds and to find examples
of manifolds $M$ with non-trivial 
$H^*_{\rho_{\lambda \omega}}(M,{\mathbb C}), \lambda \omega \ne 0$.
One of the first observations that was made in this direction:
{\it for a nilmanifold} $G/\Gamma$ {\it the cohomology
$H^*_{\lambda \omega}(G/\Gamma, {\mathbb C})$ coincides with
the cohomology $H^*_{\lambda \omega}(\mathfrak{g})$ 
associated with the
one-dimensional representation of the Lie algebra} 
$\rho_{\lambda \omega}: \mathfrak{g} \to {\mathbb C},
\rho_{\lambda \omega}(\xi) = \lambda \omega(\xi)$ and
hence $H^*_{\lambda \omega}(\mathfrak{g})=0$ by Dixmier's
theorem ~\cite{D}(Corollary \ref{dixm}).

Applying Hattori's theorem ~\cite{H}  one can observe that the isomorphism
$H^*_{\lambda \omega}(G/\Gamma, {\mathbb C})
\cong H^*_{\lambda \omega}(\mathfrak{g})$ still holds on for
compact solvmanifolds $G/\Gamma$ with completely solvable Lie
group $G$.
A kind of minimal model of a solvable Lie algebra $\mathfrak{g}$,
the free $d$-algebra
$(\Lambda^*(\omega_1, \dots, \omega_n), d)$
with differential $d$ 
$$
d \omega_i=0, i=1,\dots, k; \quad
d \omega_{j} = \alpha_{j} \wedge \omega_{j} +
P_{j}(\omega_1,\dots, \omega_{j{-}1}), j=k{+}1, \dots, n.
$$
is considered (Lemma \ref{mlemma}).
By means of $(\Lambda^*(\omega_1, \dots, \omega_n), d)$
a spectral sequence $E_r$ that converges to the
$H^*_{\lambda \omega}(\mathfrak{g})$
is constructed. $E_r$ degenerates at the first term $E_1$ if
$-\lambda \omega$ doesn't belong to the finite set
$\Omega_{\mathfrak{g}} \subset H^1(\mathfrak{g})$.
$\Omega_{\mathfrak{g}}$ is
defined by the collection
$\{\alpha_{k{+}1}, \dots, \alpha_{n} \}$ of the closed
$1$-forms that have invariant sense as
the weights of completely reducible representation
associated to the restriction $ad|_{[\mathfrak{g},\mathfrak{g}]}$
of adjoint representation $ad: \mathfrak{g} \to \mathfrak{g}$
(Theorem \ref{solv_alg_1}).

The 
main result of this article (Theorem \ref{main}):
{ \it the cohomology with local coefficients
$H^*_{\lambda \omega}(G/\Gamma, {\mathbb C})$
of a compact solvmanifold $G/\Gamma$, where $G$ is completely solvable Lie
group is non-trivial 
if and only if $-\lambda [\omega] \in
\tilde \Omega_{\mathfrak{g}}$, where $ \tilde \Omega_{\mathfrak{g}}$
is the finite subset in $H^1(G/\Gamma, {\mathbb C})$
well defined in terms of $\mathfrak{g}$.}

The author is grateful to L. Alania for helpful discussions and
attention to this work.

\section{Deformed deRham complex and Morse-Novikov theory}
Let us consider a closed compact $C^{\infty}$-manifold $M$ and
its deRham complex $(\Lambda^*(M),d)$ of differential forms. 
Let $\omega$ be a closed $1$-form on $M$ and $\lambda \in {\mathbb R}$.
Now one can define a new algebraic complex
$(\Lambda^*(M), d_{\lambda \omega})$ with a
deformed differential
$$d_{\lambda \omega}=d+\lambda \omega : \Lambda^*(M) \to \Lambda^*(M)$$
i.e. for any form $a \in \Lambda^*(M)$:
$$d_{\lambda \omega}(a)=da+\lambda \omega \wedge a.$$
Now taking $\lambda \in {\mathbb C}$ and considering the complexification 
$\Lambda^*(M) \otimes {\mathbb C}$ we come to the following important

\begin{lemma}[~\cite{N3}, ~\cite{Pa}]
1) For a closed $1$-form $d\omega=0$ the cohomology
$H^*_{\lambda \omega}(M,{\mathbb C})$
of the algebraic complex
$(\Lambda^*(M) \otimes {\mathbb C}, d_{\lambda \omega})$
coincides with the cohomology $H^*_{\rho_{\lambda \omega}}(M,{\mathbb C})$ 
with coefficients in local system of groups ${\mathbb C}$ defined by the representation
$\rho_{\lambda \omega}: \pi_1(M) \to {\mathbb C}^*$ of fundamental
group defined by the formula
$$
\rho_{\lambda \omega}([\gamma])=\exp{\int_{\gamma} \lambda \omega},
\quad [\gamma] \in \pi_1(M),
$$

2) For any pair $\omega, \omega'$ of $1$-forms such that
$\omega- \omega'=d \phi, \phi \in C^{\infty} (M)$ the cohomology
$H^*_{\lambda \omega}(M,{\mathbb C})$ and
$H^*_{\lambda \omega'}(M,{\mathbb C})$
are isomorphic to each other. This
isomorphism can be given by the gauge transformation
$$
a \to e^{\lambda \phi}a; \quad d \to e^{\lambda \phi} d e^{-\lambda \phi}=
d+ \lambda d \phi \wedge
$$
\end{lemma}

We denote corresponding Betti numbers by
$b_p(\lambda, \omega)$, where
$b_p(\lambda, \omega)=\dim H^*_{\rho_{\lambda \omega}}(M,{\mathbb C})$.

\begin{remark}
The representation $\rho_{\lambda \omega}: \pi_1(M) \to
{\mathbb C}^*$ 
defines a local system of groups ${\mathbb C}$
on the manifold $M$ in the sense of Steenrod (see ~\cite{R} for details).
\end{remark}

The cohomology $H^*_{\lambda \omega}(M,{\mathbb C})$ naturally arises
in the Morse-Novikov theory:
we assume now that $\omega$ is a Morse $1$-form, i.e., in a neighbourhood
of any point $\omega=df$, where $f$ is a Morse function.
The zeros of $\omega$ are isolated, and one can
define the index of each zero. The number of zeros of $\omega$
of index $p$ is denoted by $m_p(\omega)$.

\begin{theorem}[A. Pazhitnov, \cite{Pa}]
For sufficiently large real numbers $\lambda$,
$$ m_p(\omega) \ge b_p(\lambda, \omega) \; \forall p.$$
\end{theorem}

\begin{theorem}[A. Pazhitnov, \cite{Pa}]
Assume that all the periods of the form $\omega$ are commensurable.
If $Re \lambda$ is sufficiently large, then 
$b_p(\lambda, \omega)=b_p(\omega)$, where $b_p(\omega)$ is a Novikov number.
\end{theorem}

\section{Dixmier's exact sequence of Lie algebra cohomology}

Let $\mathfrak{g}$ be a $n$-dimensional Lie algebra.
The dual of the Lie bracket
$[,]: \Lambda^2 (\mathfrak{g}) \to \mathfrak{g}$ gives a linear mapping
$d_1: \mathfrak{g}^* \to \Lambda^2 (\mathfrak{g}^*)$
which extends in a standard way to a differential $d$
of a cochain complex of the Lie algebra
$\mathfrak{g}$:
$$
\begin{CD}
{\mathbb K} @>{d_0{=}0}>>
\mathfrak{g}^* @>{d_1}>> \Lambda^2 (\mathfrak{g}^*) @>{d_2}>>
\Lambda^3 (\mathfrak{g}^*) @>{d_3}>>
\dots \end{CD}
$$

$$
d(\rho \wedge \eta)=d\rho \wedge \eta+(-1)^{deg\rho} \rho \wedge d\eta,
\; \forall \rho, \eta \in \Lambda^{*} (\mathfrak{g}^*).
$$

Vanishing of the $d^2$ corresponds to the Jacobi identity.

For $d:\Lambda^q(\mathfrak{g}^*) \to \Lambda^{q{+}1}(\mathfrak{g}^*)$
and $f \in \Lambda^q(\mathfrak{g}^*)$ the following formula holds on:

\begin{equation}
\begin{split}
d f(X_1, \dots, X_{q{+}1})= \\
\sum_{1{\le}i{<}j{\le}q{+}1}({-}1)^{i{+}j{-}1}
f([X_i,X_j],X_1, \dots, \hat X_i, \dots, \hat X_j, \dots, X_{q{+}1}).
\end{split}
\end{equation}

A cohomology of this complex is called
the cohomology (with trivial coefficients) of the Lie algebra
$\mathfrak{g}$ and is denoted by $H^*(\mathfrak{g})$.

From the definition 
it follows that
$H^1(\mathfrak{g})$ is the dual space to
$\mathfrak{g}/[\mathfrak{g},\mathfrak{g}]$ and so

1) $b^1(\mathfrak{g})=\dim H^1(\mathfrak{g}) \ge 2$ for
a nilpotent Lie algebra $\mathfrak{g}$, 

2) $b^1(\mathfrak{g}) \ge 1$ for
a solvable Lie algebra $\mathfrak{g}$,

3) $b^1(\mathfrak{g})=0$
for a semi-simple Lie algebra $\mathfrak{g}$.

Now we take a Lie algebra $\mathfrak{g}$
over a field ${\mathbb K}$ with a non-trivial $H^1(\mathfrak{g})$.
Let $\omega \in \mathfrak{g}^*, \omega \ne 0, d \omega =0$ and
$\lambda \in {\mathbb K}$. One can define

1) a new deformed differential $d_{\lambda \omega}$
in $\Lambda^{*} (\mathfrak{g}^*)$
by the formula
$$d_{\lambda \omega}(a)=da+\lambda \omega \wedge a.$$

2) an one-dimensional representation
$$\rho_{\lambda \omega}: \mathfrak{g} \to {\mathbb K}, \:
\rho_{\lambda \omega}(\xi)= \lambda \omega(\xi), \xi \in \mathfrak{g}.$$

Now we recall the definition of the Lie algebra cohomology associated
with a representation. Let $\mathfrak{g}$ be a Lie algebra and
$\rho: \mathfrak{g} \to \mathfrak{gl}(V)$ its linear representation.
We denote by $C^q(\mathfrak{g},V)$
the space of $q$-linear alternating mappings of $\mathfrak{g}$ into
$V$. Then one can consider an algebraic complex:

$$
\begin{CD}
V=C^0(\mathfrak{g}, V) @>{d_0}>>
C^1(\mathfrak{g}, V) @>{d_1}>> C^2(\mathfrak{g}, V) @>{d_2}>>
\dots
\end{CD}
$$
where the differential $d_q$ is defined by:

\begin{equation}
\begin{split}
(d_q f)(X_1, \dots, X_{q{+}1})=
\sum_{i{=}1}^{q{+}1}(-1)^{i{+}1}
\rho(X_i)(f(X_1, \dots, \hat X_i, \dots, X_{q{+}1}))+\\
+ \sum_{1{\le}i{<}j{\le}q{+}1}(-1)^{i{+}j{-}1}
f([X_i,X_j],X_1, \dots, \hat X_i, \dots, \hat X_j, \dots, X_{q{+}1}).
\end{split}
\end{equation}

The cohomology of the complex $(C^*(\mathfrak{g}, V), d)$ is called
the cohomology of the Lie algebra $\mathfrak{g}$
associated to the representation $\rho: \mathfrak{g} \to V$.

\begin{proposition}
Let $\mathfrak{g}$ be a Lie algebra and $\omega$
is a closed $1$-form.
Then the complex $(\Lambda^*(\mathfrak{g^*}), d_{\lambda \omega})$
coincides with the cochain complex of the Lie algebra
$\mathfrak{g}$ associated with one-dimensional representation
$\rho_{\lambda \omega}: \mathfrak{g} \to {\mathbb K}$,
where $\rho_{\lambda \omega}(\xi)=\lambda \omega(\xi), \xi \in \mathfrak{g}$.
\end{proposition}

The proof follows from
$$(\lambda \omega \wedge a) (X_1, \dots, X_{q{+}1})=
\sum_{i{=}1}^{q{+}1}(-1)^{i{+}1}
\lambda \omega(X_i)(a(X_1, \dots, \hat X_i, \dots, X_{q{+}1})).
$$

One can deduce that
$H^{0}_{\lambda \omega}(\mathfrak{g})=0$ for a non-trivial
$\lambda \omega$, as well as $H^{n}_{\lambda \omega}(\mathfrak{g})=0$
for an unimodular $n$-dimensional Lie algebra $\mathfrak{g}$.

\begin{remark}
The cohomology $H^{*}_{\lambda \omega}(\mathfrak{g})$
coincides with the Lie algebra cohomology $H^{*}(\mathfrak{g})$ with trivial coefficients if
$\lambda=0$. If $\lambda \ne 0$
the deformed differential $d_{\lambda \omega}$
is not compatible with the exterior product $\wedge$
in $\Lambda^*(\mathfrak{g})$
$$d_{\lambda \omega}(a \wedge b)=d(a \wedge b)+
\lambda \omega \wedge a \wedge b
\ne d_{\lambda \omega}(a) \wedge b +
(-1)^{deg a}a \wedge d_{\lambda \omega}(b) $$
and
the 
cohomology $H^{*}_{\lambda \omega}(\mathfrak{g})$ has no
natural multiplicative structure
and therefore no Poincare
duality in the case of unimodular Lie algebra $\mathfrak{g}$.
The corresponding Euler characteristic
$\chi_{\lambda \omega}(\mathfrak{g})$
is still equal to zero.
\end{remark}

Let $\omega \in \mathfrak{g}^*, \omega \ne 0, d \omega = 0$.
Then $\mathfrak{b}_{\omega}=\{x \in \mathfrak{g}, \omega(x)=0 \}$ is an
ideal of codimension $1$ in $\mathfrak{g}$.
One can choose an element
$X\in \mathfrak{g}, \omega(X)=1$. 

\begin{theorem}[Dixmier ~\cite{D}]
There exists a long exact sequence of Lie algebra cohomology:
\begin{equation}
\dots 
\stackrel{adX^*_{i{-}1}{+}\lambda Id}{\longrightarrow}
H^{i{-}1}(\mathfrak{b}_{\omega})
\stackrel{\omega \wedge}{\longrightarrow}
H^{i}_{\lambda \omega}(\mathfrak{g})
\stackrel{r_i}{\rightarrow} H^{i}(\mathfrak{b}_{\omega})
\stackrel{adX^*_i{+}\lambda Id}{\longrightarrow} H^{i}(\mathfrak{b}_{\omega})
\to \dots
\end{equation}
where

1) the homomorphism
$r_i: H^{i}(\mathfrak{g}) \to H^{i}(\mathfrak{b})$
is the restriction homomorphism;

2) $\omega \wedge: H^{*{-}1}(\mathfrak{b}_{\omega})
\to H^{*}(\mathfrak{g})$
is induced by the multiplication
$\omega \wedge: \Lambda^{*{-}1}(\mathfrak{b}_{\omega}^*)
\to \Lambda^*(\mathfrak{g}^*)$;

3) the homomorphisms
$adX^*_i: H^{i}(\mathfrak{b}_{\omega}) \to H^{i}(\mathfrak{b}_{\omega})$
are induced by the derivation $adX^*_i$ of degree 
zero of $\Lambda^*(\mathfrak{b}^*_{\omega})$ ($adX^*(a\wedge b)=
adX^*a\wedge b+ a\wedge adX^*b, \; \forall a, b \in
\Lambda^*(\mathfrak{b}^*_{\omega})$) that continues a dual mapping
$adX^* : \mathfrak{b}^*_{\omega} \to \mathfrak{b}^*_{\omega}$ to the
$adX : \mathfrak{b}_{\omega} \to \mathfrak{b}_{\omega}$ operator.
The derivation $adX^*$ commutes
with $d$ and corresponding mapping in $H^{*}(\mathfrak{b}_{\omega})$
we denote
by the same symbol $adX^*$. Id is the identity operator.
\end{theorem}

Each form $f \in \Lambda^*(\mathfrak{g})$
can be decomposed
as $f=\omega \wedge f' + f''$, where $ f' \in \Lambda^{*{-}1}(\mathfrak{b}^*_{\omega})$ and 
$ f'' \in \Lambda^{*}(\mathfrak{b}^*_{\omega})$.
And one can write out a short exact sequence of algebraic complexes
$$
0 \longrightarrow \Lambda^{*{-}1}(\mathfrak{b}^*_{\omega})
\stackrel{\omega \wedge}{\longrightarrow} \Lambda^{*}(\mathfrak{g}^*)
\longrightarrow \Lambda^{*}(\mathfrak{b}^*_{\omega})
\longrightarrow 0
$$
where
$\Lambda^{*}(\mathfrak{b}^*_{\omega})$ has
the standard differential $d$,
$\Lambda^{*}(\mathfrak{g}^*)$ has the deformed differential
$d_{\lambda \omega}$
and
$\Lambda^{*{-}1}(\mathfrak{b}^*_{\omega})$
is taken with the differential $-d$ as
$d(\omega \wedge c)=d_{\lambda \omega}(\omega \wedge c)=
-\omega \wedge dc$. 

The short exact sequence of algebraic complexes gives us
the long exact sequence of cohomology. Everything is clear with the homomorphisms
$r_i$ and $\omega \wedge$. As for
the homomorphism
$H^{q}(\mathfrak{b}_{\omega}) \to H^{q}(\mathfrak{b}_{\omega})$
we can remark the following.

First of all let us introduce a new mapping
$$
\Lambda^{*}(\mathfrak{g}^*)
\to \Lambda^{*-1}(\mathfrak{b}^*_{\omega})
,
f \in \Lambda^{*}(\mathfrak{g}^*)\to f_X \in
\Lambda^{*-1}(\mathfrak{b}_{\omega}^*),$$ where
$f_X(X_1,\dots,X_q)=f(X,X_1,\dots,X_q)$. So $f_X= f'$ if $f=\omega \wedge f' +f''$.

Then an obvious formula holds on:

\begin{multline}
(d f)_X(X_1, \dots, X_{q{+}1})=
\sum_{1{\le}j{\le}q{+}1}({-}1)^{j}
f(adX(X_j), X_1,\dots, \hat X_j, \dots, X_{q{+}1})+\\
+\sum_{1{\le}i{<}j{\le}q{+}1}({-}1)^{i{+}j}
f([X_i,X_j],X,X_1, \dots, \hat X_i, \dots, \hat X_j, \dots, X_{q{+}1})=\\
=(adX^*_q(f)+d(f_X))(X_1, \dots, X_{q{+}1}).
\end{multline}

Hence the homomorphism $H^{q}(\mathfrak{b}_{\omega}) \to H^{q}(\mathfrak{b}_{\omega})$,
$[f] \to [f' + \lambda f]$ of long exact sequence in cohomology coincides with homorphism induced by 
$adX^* +\lambda Id$.

\begin{corollary}
\label{cor_D}
Let $\mathfrak{g}$ be a $n$-dimensional Lie algebra and
$\omega \in \mathfrak{g}^*, \omega \ne 0, d \omega=0$
and $Spec^k (\omega)$  be the set of eigenvalues of operator
$adX^*_k: H^{k}(\mathfrak{b}_{\omega}) \to H^{k}(\mathfrak{b}_{\omega})$,
then

1) the cohomology $H^{*}_{\lambda \omega}(\mathfrak{g})$
is non-trivial if and only if 
$$-\lambda \in \cup_{k{=}1}^n Spec^k (\omega);$$

2) the $i$-th Betti number
$b^i_{\lambda \omega}(\mathfrak{g})=dim H^i_{\lambda \omega}(\mathfrak{g})$
can be expressed in a following way:
\begin{equation}
b^i_{\lambda \omega}(\mathfrak{g})=
k^i_{\lambda \omega}+k^{i{-}1}_{\lambda \omega},
\end{equation}
where by $k^i_{\lambda \omega}$ we denote the dimension of the
kernel of $adX^*_i+\lambda Id: H^{i}(\mathfrak{b}_{\omega}) \to
H^{i}(\mathfrak{b}_{\omega})$.

\end{corollary}

\begin{example}
Let $\mathfrak{g}$ be a Lie algebra defined by the basis
$X, e_1, e_2, \dots, e_n$
and commutating relations (trivial ones are omitted):
$$
[X,e_1]=e_1, [X,e_2]=e_2, \dots, [X,e_n]=e_n.
$$
Thus $\mathfrak{g}$ is a semidirect sum of ${\mathbb K}$
and abelian ${\mathbb K}^n$ defined by the operator $adX$
with the identity matrix $E$ in the basis
$e_1, e_2, \dots, e_n$ of ${\mathbb K}^n$.
We take $\omega, \omega_1, \dots, \omega_n$ as the corresponding dual basis in $\mathfrak{g}$. In particular 
$$\omega(X)=1, \omega(e_i)=0, i=1, \dots, n.$$ $d \omega=0$
and $adX^*(\omega_{i_1}\wedge\dots\wedge \omega_{i_p})=p \omega_{i_1}\wedge \dots \wedge
\omega_{i_p}$. Hence $Spec^p(\omega)=\{p\}$ and for $\lambda=p$
we have only two
non-trivial Betti numbers: $$b^{p{-}1}_{p \omega}(\mathfrak{g})=
b^{p}_{p \omega}(\mathfrak{g})=\binom{n}{p}.$$
\end{example}

\begin{corollary}[Dixmier, ~\cite{D}]
\label{dixm}
Let $\mathfrak{g}$ be a nilpotent Lie algebra.
The cohomology $H^{*}_{\lambda \omega}(\mathfrak{g})$
is trivial for all non-trivial $\lambda\omega$ and coincides
with the Lie algebra cohomology $H^{*}(\mathfrak{g})$ in trivial $\lambda\omega=0$ case.
\end{corollary}

The operator $adX$ is nilpotent
and therefore the same is
$adX^*_i: H^{i}(\mathfrak{b}_{\omega}) \to H^{i}(\mathfrak{b}_{\omega})$.
Hence
$adX^*_i{+}\lambda Id$ is non-degenerate operator for all $\lambda \neq 0$.
We obtain the proof of Theorem 1 from ~\cite{Al} as
the corollary of Dixmier's theorem ~\cite{D} for cohomology
of nilpotent Lie algebras.

\begin{remark}
We represented in this article only a version of Dixmier's
exact sequence adapted to our special case of
$1$-dimensional Lie algebra representation (see ~\cite{D}
for all details), the last thing that we want to recall
is Dixmier's estimate for Betti numbers $\dim H^p(\mathfrak{g})$ of a nilpotent
Lie algebra $\mathfrak{g}$.
\end{remark}

\begin{corollary}[~\cite{D}]
Let $\mathfrak{g}$ be a nilpotent Lie algebra. Then
$$\dim H^p(\mathfrak{g}) \ge 2, \; p=1,\dots, n{-}1.$$
\end{corollary}

It follows from the Corollary \ref{cor_D}: we have $\lambda \omega=0$
and $\dim H^p(\mathfrak{g})=k^p+k^{p{-}1}$, where
$k^p \ge 1, k^{p{-}1} \ge 1$ are the dimensions of the kernels of nilpotent
operators $adX^*_p$, $adX^*_{p{-}1}$ in the spaces
$H^p(\mathfrak{b}_{\omega})$, $H^{p{-}1}(\mathfrak{b}_{\omega})$
with $\dim H^p(\mathfrak{b}_{\omega}) \ge 2$,
$\dim H^{p{-}1}(\mathfrak{b}_{\omega}) \ge 2$
by inductive assumption.

\section{Cohomology of solvable Lie algebras}

\begin{definition}
A real solvable Lie algebra $\mathfrak{g}$ is called completely solvable
if $ad(X): \mathfrak{g} \to \mathfrak{g}$
has only real eigenvalues $\forall X \in \mathfrak{g}$.
\end{definition}

\begin{lemma}
\label{mlemma}
Let $\mathfrak{g}$ be a $n$-dimensional solvable  over ${\mathbb C}$
(or real completely solvable Lie algebra) and
$b^1(\mathfrak{g})=\dim H^1(\mathfrak{g})=k$.
Then exists a basis $\omega_1, \dots, \omega_n$
in $\mathfrak{g}^*$ such that
\begin{multline}
\label{solvemodel}
d \omega_1= \dots = d \omega_k=0,\\
d \omega_{k{+}1} = \alpha_{k{+}1} \wedge \omega_{k{+}1} +
P_{k{+}1}(\omega_1,\dots, \omega_{k}), \\
\dots, \\
d \omega_{n} = \alpha_{n} \wedge \omega_{n} +
P_{n}(\omega_1,\dots, \omega_{n{-}1}),
\end{multline}
where 
$\alpha_{k{+}1}, \dots, \alpha_{n}$ are closed $1$-forms, that are
the weights of completely reducible representation
associated to $ad|_{[\mathfrak{g}, \mathfrak{g}]}$
and $P_{i}(\omega_1,\dots, \omega_{i{-}1}) \in
\Lambda^2(\omega_1,\dots, \omega_{i{-}1})$.
\end{lemma}

For the proof we apply Lie's theorem to
the adjoint representation $ad$ restricted 
to the commutant $[\mathfrak{g}, \mathfrak{g}]$.
Namely we can choose a basis $e_{k{+}1}, \dots, e_{n}$ such that
the subspaces $Span(e_{i}, \dots, e_{n}), i=k{+}1, \dots, n$
are invariant with respect to the representation
$ad|_{[\mathfrak{g}, \mathfrak{g}]}$. Then we add $e_{1}, \dots, e_{k}$
in a way that $e_{1}, \dots, e_{n}$ form the basis of $\mathfrak{g}$.
For the dual forms $\omega_{1}, \dots, \omega_{n}$
in $\mathfrak{g}^*$ we have formulas ~(\ref{solvemodel}). 

\begin{remark}
One can consider the free $d$-algebra
$(\Lambda^*(\omega_1, \dots, \omega_n), d)$ as a kind of minimal model
of the cohaine complex $\Lambda^*(\mathfrak{g}^*)$ because the mapping
$\Lambda^*(\omega_1, \dots, \omega_n) \to \Lambda^*(\mathfrak{g}^*)$
induces the isomorphism in cohomology.
\end{remark}

Now we start with a solvable 
Lie algebra over ${\mathbb C}$ (or real completely solvable)
and take a basis $\omega_{1}, \dots, \omega_{n}$
constructed in Lemma ~\ref{mlemma}. 
Let us denote by $\Lambda^*$ the exterior
subalgebra in $\Lambda^*(\omega_1, \dots, \omega_n)$
generated by $\omega_1, \dots, \omega_k$.
One can define
a filtration $F$ of the cochain complex $\Lambda^*(\mathfrak{g}^*)$:
$$
0 \subset \Lambda^* \subset F^{k{+}1} 
\subset F^{k{+}2} \subset F^{k{+}1,k{+}2} \subset F^{k{+}3}
\subset F^{k{+}1,k{+}3} \subset 
\dots
\subset F^{k,\dots,n}=\Lambda^*(\mathfrak{g}^*)
$$
where the system of subspaces
$\{ F^{j_1,\dots,j_p}, k < j_1 < \dots <j_p \le n \}$
is defined by the following conditions:

A subspace $F^{j_1,\dots,j_p}$ is spanned by monomials
$a=\omega_{l_1} \wedge \dots \wedge \omega_{l_q}$ with  $l_1 < \dots <l_q$
such that 

$l_q < j_p$; or $l_q = j_p,
l_{q{-}1}<j_{p{-}1}$; 

or $l_q = j_p, l_{q{-}1} = j_{p{-}1}, l_{q{-}2}<j_{p{-}2}$; or $\dots$;

or $l_q = j_p, \dots, l_{q{-}p{+}1} = j_{1}, l_{q{-}p} \le k$.

Thus for example:
\begin{multline}
F^{k{+}1}=\Lambda^* \oplus \Lambda^* {\wedge} \omega_{k{+}1}, 
F^{k{+}2}=\Lambda^* \oplus \Lambda^* {\wedge} \omega_{k{+}1}
\oplus \Lambda^* {\wedge} \omega_{k{+}2}, \\
F^{k{+}1,k{+}2}=\Lambda^* \oplus \Lambda^* {\wedge} \omega_{k{+}1}
\oplus \Lambda^* {\wedge} \omega_{k{+}2}
\oplus \Lambda^* {\wedge} \omega_{k{+}1} {\wedge} \omega_{k{+}2},\\
F^{k{+}3}=\Lambda^* \oplus \Lambda^* {\wedge} \omega_{k{+}1}
\oplus \Lambda^* {\wedge} \omega_{k{+}2}
\oplus \Lambda^* {\wedge} \omega_{k{+}1} {\wedge} \omega_{k{+}2}
\oplus \Lambda^* {\wedge} \omega_{k{+}3}, \dots
\end{multline}

The filtration $F$ is compatible with differential
$d+\lambda \omega$ and one can consider the corresponding
spectral sequence $E_r$. To obtain its first term $E_1$ one
have to calculate the cohomology of complexes
$(\Lambda^* \wedge \omega_{j_1} {\wedge} \dots {\wedge} \omega_{j_p}
, d_0)$:
$$d_0(\tilde a\wedge \omega_{j_1} {\wedge} \dots {\wedge} \omega_{j_p})
=(\alpha_{j_1}{+}\dots{+}\alpha_{j_p}{+}\lambda\omega)
\tilde a \wedge \omega_{j_1} {\wedge} \dots {\wedge} \omega_{j_p}.$$
The cohomology
$H^*(\Lambda^* \wedge \omega_{j_1} {\wedge} \dots
{\wedge} \omega_{j_p}, d_0)$
coincides with the
cohomology of $(\Lambda^*, \tilde d)$ where differential $\tilde d$ acts
as exterior multiplication by $1$-form
$\alpha_{j_1}{+}\dots{+}\alpha_{j_p}{+}\lambda\omega$.
Hence
$H^*(\Lambda^* \wedge \omega_{j_1} {\wedge} \dots {\wedge} \omega_{j_p}, d_0)$
is trivial if
$\alpha_{j_1}{+}\dots{+}\alpha_{j_p}{+}\lambda\omega \ne 0$.
So taking $\lambda \omega$ such that 
$$\alpha_{j_1}{+}\dots{+}\alpha_{j_p}{+}\lambda\omega \ne 0 \; \forall
j_1{<}\dots{<}j_p, p=k{+}1, \dots, n$$
one can see that the spectral sequence $E_r$
degenerates at the first term $E_1$.
Taking the complexification of
real solvable Lie algebra $\mathfrak{g}$ we 
come to the following
\begin{theorem}
\label{solv_alg_1}
Let $\mathfrak{g}$ be a solvable Lie algebra, $\dim \mathfrak{g}=n$
and $\{ \alpha_{k{+}1}, \dots, \alpha_{n}  \}$
is the collection of the weights of completely reducible representation
associated to $ad|_{[\mathfrak{g}, \mathfrak{g}]}$.
Let $\Omega_{\mathfrak{g}}$  denote
the set of all $p$-sums 
$\alpha_{i_1}{+} \dots{+} \alpha_{i_p}$,
$k{+}1 \le i_1{<} \dots {<} i_p \le n$, $p=1, \dots, n$. Then
exists a spectral sequence $E_r$ that converges 
to the cohomology $H^*_{\lambda \omega}(\mathfrak{g})$ and
its first term $E_1$ degenerates
if $-\lambda \omega \notin \Omega_{\mathfrak{g}}$.
\end{theorem}

The set $\Omega_{\mathfrak{g}}$ is defined by $\lambda \omega$ such that
the term $E_1$ is non-trivial, but generally $E_1$ doesn't coincide with
$E_{\infty}$. So
we have to introduce
$\tilde \Omega_{\mathfrak{g}} \subset \Omega_{\mathfrak{g}}$
such that 
$E_{\infty} \ne 0$ if and only if $ -\lambda \omega \in \{0\} \cup \tilde \Omega_{\mathfrak{g}}$.

\begin{corollary}
\label{solv_alg_2}
Let $\mathfrak{g}$ be a solvable Lie algebra.
Then $H^*_{\rho_{\lambda \omega}}(\mathfrak{g})$ is non-trivial
if and only if
$-\lambda \omega  \in \{ 0 \} \cup \tilde \Omega_{\mathfrak{g}}$ --
the finite subset in $H^1(\mathfrak{g})$.
\end{corollary}

\section{Cohomology of solvmanifolds.}

\begin{definition}  A solvmanifold (nilmanifold) $M$
is a compact homogeneous space
of the form
$G/\Gamma,$ where $G$ is a simply connected solvable (nilpotent) Lie group
and $\Gamma$ is a lattice in $G$.
\end{definition}

Let $\mathfrak{g}$ denote a Lie algebra of $G$. Recall that
$G$ is solvable if and only if $\mathfrak{g}$ is solvable Lie
algebra, the last condition is equivalent to nilpotency of
derived algebra $[\mathfrak{g},\mathfrak{g}]$.

We start with examples of nilmanifolds.

\begin{example}
A $n-$dimensional torus $T^n={\mathbb R}^n/{\mathbb Z}^n$.
\end{example}
\begin{example}
The Heisenberg manifold $M_3={\mathcal H}_3/\Gamma_3,$ where
${\mathcal H}_3$ is the group of all matrices of the form
$$
   \left( \begin{array}{lcr}
   1 & x & z\\
   0 & 1 & y\\
   0 & 0 & 1\\
   \end{array} \right) , ~~~ x,y,z \in {\mathbb R},
$$
and a lattice $\Gamma_3$ is a subgroup of matrices with $x,y,z
\in {\mathbb Z}$.
\end{example}

\begin{theorem}[A.I. Malcev \cite{Mal}]
Let $G$ be a simply connected nilpotent Lie group with a tangent Lie
algebra $\mathfrak{g}$.
Then $G$ has a co-compact lattice $\Gamma$ (i.e. $G/\Gamma$ is a compact
space) if and only if there exists a 
basis $e_1, e_2, \dots, e_n$ in $\mathfrak{g}$
such that the constants $\{c_{ij}^{k} \}$ of Lie structure
$[e_i,e_j]= c_{ij}^{k}e_k$
are all rational numbers.
\end{theorem}

This theorem gives us a practical tool for construction of
nilmanifolds: let $\mathfrak{g}$ be a nilpotent Lie algebra
defined by its basis
$e_1, e_2, \dots, e_n$
and commutating relations $[e_i,e_j]= c_{ij}^{k}e_k$, where
all numbers $c_{ij}^{k} \in {\mathbb Q}$.
Now one can define
a group structure $*$ in the vector space $\mathfrak{g}$
using the Campbell-Hausdorff formula.
The nilpotent group $G=(\mathfrak{g},*)$ 
has a co-compact lattice $\Gamma$ (a subgroup generated
by basic elements $e_1, e_2, \dots, e_n$) and one can consider
corresponding nilmanifold $G/\Gamma$.

\begin{example}
Let ${\mathcal V}_n$ be a nilpotent Lie algebra
with a basis $e_1, e_2, \dots, e_n$ and a Lie bracket:
$$[e_i,e_j]= \left\{\begin{array}{r}
   (j-i)e_{i{+}j},  i+j \le n;\\
   0, i+j > n.\\
   \end{array} \right. $$
Cohomology of the corresponding family of nilmanifolds
$M_n$ was studied in ~\cite{Al} and ~\cite{Mill}.
\end{example}

The situation with non-nilpotent solvable Lie groups is
much more difficult: the crucial point is
the problem of existence of cocompact lattice (see ~\cite{R}
for details).
For example, if a solvable Lie group $G$ admits a cocompact lattice $\Gamma$
then the corresponding Lie algebra $\mathfrak{g}$ is unimodular,
hence $\alpha_{k{+}1}+\dots+\alpha_{n}=0$. The condition
of unimodularity of $\mathfrak{g}$ is not sufficient.
See ~\cite{R} for general information in this domain.

\begin{example}
\label{3-dim}
Let us consider a semidirect product
$G_0={\mathbb R} \ltimes {\mathbb R}^2$ where ${\mathbb R}$ acts
on ${\mathbb R}^2$ via
$$
t \to \phi(t)=\begin{pmatrix}a^t&0\\
0 & a^{{-}t} \end{pmatrix},
$$
where $a+a^{-1}=n \in {\mathbb N}, a \ne 1, a > 0$.
Then
$$
\phi(1)= C^{-1}\begin{pmatrix}0&1\\
-1 & n \end{pmatrix} C,
$$
for some matrix $C \in GL(2, {\mathbb R})$. Then exists a lattice
$L \subset {\mathbb R}^2$ invariant with respect to $\phi(1)$.
${\mathbb Z}$ acts on $L$ via $\phi(1)^m$ and we define
a lattice $\Gamma \subset G$ as ${\mathbb Z} \ltimes_{\phi(1)} L$.
The lattices corresponding to different values of $n$
are, generally speaking, non-isomorphic. So corresponding solvmanifolds
are non diffeomorphic to. But we study the cohomology $H^*_{\rho}(G/\Gamma, {\mathbb C})$ over 
${\mathbb C}$ and as we  will see it doesn't depend on the choice of $\Gamma \subset G$.
\end{example}

The solvmanifold from the previous example is a fibre bundle over $S^1$ with  $T^2$ as fibre.
It can be generalized by the following  

\begin{theorem}[G.D. Mostow ~\cite{Mos1}]
Any compact solvmanifold is a bundle with toroid as base space
and nilmanifold as fibre.
\end{theorem}

\begin{definition}
A solvable Lie group $G$ is called completely solvable if
its tangent Lie algebra $\mathfrak{g}$ is completely solvable.
\end{definition}

One can identify deRham complex $\Lambda^*(G/\Gamma)$
with subcomplex
$\Lambda^*_{\Gamma}(G) \subset \Lambda^*(G)$
of left-invariant forms on $G$ with respect to the action of the lattice
$\Gamma$. $\Lambda^*_{\Gamma}(G)$ containes the
subcomplex $\Lambda^*_{G}(G)$ of left-invariant forms 
with respect to the whole action of $G$.
$\Lambda^*_{G}(G)$ is naturally isomorphic
to the Lie algebra cochain complex
$\Lambda^*(\mathfrak{g})$.
Let us consider the corresponding inclusion
$$\psi: \Lambda^*(\mathfrak{g}) \to 
\Lambda^*(G/\Gamma).$$

\begin{theorem}[Hattori ~\cite{H}]
Let $G/\Gamma$ be a compact solvmanifold, where
$G$ is a completely solvable Lie group, then
the inclusion $\psi: \Lambda^*(\mathfrak{g}) \to 
\Lambda^*(G/\Gamma)$ induces the isomorphism 
$\psi^*: H^*(\mathfrak{g}) \to H^*(G/\Gamma, {\mathbb R})$
in cohomology.
\end{theorem}

\begin{remark}
In fact Hattori's theorem is the generalization of the theorem
proved by Nomizu ~\cite{Nz} for nilmanifolds.
For an arbitrary solvmanifold $G/\Gamma$ the mapping $\psi^*$
is not isomorphism but it is
an inclusion (see ~\cite{R}).
\end{remark}

So every class $[\omega] \in H^1(G/\Gamma, {\mathbb R})$ can be
represented by the
left-invariant (with respect to the action of $G$)
$1$-form $\omega$. By means of $\omega$ one can define
a one-dimenisional representation
$\rho_{\lambda \omega}: G \to {\mathbb C}^*$:
$$ \rho_{\lambda \omega}(g)=\exp \int_{\gamma(e,g)} \lambda \omega,$$
where $\gamma(e,g)$ is a path connecting the identity $e$ with
$g \in G$ (let us recall that $G$ is a symply connected).
As $\omega$ is the left invariant $1$-form then
$$ \int_{\gamma(e,g_1 g_2)} \lambda \omega=
\int_{\gamma(e,g_1)} \lambda \omega+
\int_{\gamma(g_1,g_1 g_2)} \lambda \omega=
\int_{\gamma(e,g_1)} \lambda \omega+
\int_{g_1^{{-}1}\gamma(e,g_2)} \lambda \omega
$$
holds on and
$\rho_{\lambda \omega}(g_1g_2)=\rho_{\lambda \omega}(g_1)
\rho_{\lambda \omega}(g_2)$. $\rho_{\lambda \omega}$
induces the representation of corresponding Lie algebra
$\mathfrak{g}$ (we denote it by the same symbol):
$\rho_{\lambda \omega}(X)= \lambda \omega(X)$.

In this situation it's possible to make some generalizations using 
Mostow's theorem. Namely following ~\cite{R} we give

\begin{definition}
Let $G$ be a simply-connected Lie group and $\Gamma \subset G$ a lattice. Let $\rho$ be a finite dimensional
representation of $G$ on a complex vector space $F$. Let $Ad$ denote the adjoint representation of $G$ on its
Lie algebra $\mathfrak{g}$ as well as the complexification $\mathfrak{g}_{\mathbb C}$ of $\mathfrak{g}$. We will
say that the representation $\rho$ is $\Gamma$-supported if $\rho(\Gamma)$
and $\rho(G)$ have the same Zariski closure in $Aut_{\mathbb C}(F)$. The representation $\rho$ is $\Gamma$-admissible
if $\rho \oplus Ad$ (on  $F \oplus \mathfrak{g}_{\mathbb C}$) is $\Gamma$-supported.
\end{definition}

\begin{theorem}[Mostow, Theorem 7.26 in ~\cite{R}]
Let $G$ be a simply-connected solvable Lie group and
a $\Gamma \subset G$ a lattice. Let $\rho$ be a finite dimensional 
$\Gamma$-admissible representation in a complex vector space $F$.

Then the inclusion $\psi: \Lambda^*(\mathfrak{g},F) \to 
\Lambda^*(G/\Gamma,F)$ of complexes of differential forms with values in $F$ induces the isomorphism 
$\psi^*: H^*_{\rho}(\mathfrak{g},F) \to H^*_{\rho}(G/\Gamma, F)$
in cohomology where $\rho$ is used also to denote
the representation of the Lie algebra
$\rho: \mathfrak{g} \to F$ induced by $\rho: G \to Aut_{\mathbb C}(F)$.
\end{theorem}

But we'll not discuss the details of Mostow's theorem and possible generalizations
we restrict ourselves to the case of completely solvable Lie group $G$. Namely  we'll prove by means of Hattori's theorem
the following important

\begin{corollary}
\label{main_c}
Let $G/\Gamma$ be a compact solvmanifold, 
$G$ has a completely solvable Lie group and
$\omega$ is a closed $1$-form on $G/\Gamma$.
The cohomology $H^*_{\lambda \omega}(G/\Gamma, {\mathbb C})$
is isomorphic to the Lie algebra cohomology
$H^*_{\lambda \omega'}(\mathfrak{g})$ where
$\omega' \in \mathfrak{g}^*$ is
the left-invariant $1$-form that represents
the class $[\omega] \in H^1(G/\Gamma, {\mathbb R})$.
\end{corollary}

The Corollary ~\ref{main_c} together with Corollary ~\ref{solv_alg_2} gives us

\begin{theorem}
\label{main}
Let $G/\Gamma$ be a compact solvmanifold,  
$G$ is a completely solvable Lie group and
$\omega$ is a closed $1$-form on $G/\Gamma$.
The cohomology $H^*_{\lambda \omega}(G/\Gamma, {\mathbb C})$ is non-trivial if and only if
$-\lambda [\omega] \in
\{0\} \cup \tilde \Omega_{\mathfrak{g}}$ -- the finite subset in
$H^1(G/\Gamma, {\mathbb R})$ well-defined in terms of the
 corresponding Lie algebra $\mathfrak{g}$.
\end{theorem}

\begin{corollary}[~\cite{Al}]
Let $G/\Gamma$ be a compact nilmanifold.
The cohomology $H^*_{\lambda \omega}(G/\Gamma, {\mathbb C})=0$
if and only if $\lambda \omega \ne 0$.
\end{corollary}

Let us consider a $3$-dimensional solvmanifold $M=G_0/\Gamma$
defined in the Example ~\ref{3-dim}.
The corresponding Lie algebra $\mathfrak{g}_0$
is isomorphic to the Lie algebra defined by its basis $e_1,e_2,e_3$
and following non-trivial brackets:
$$[e_1,e_2]=e_2, \; [e_1,e_3]=-e_3.$$
We take a dual basis $\omega_1,\omega_2,\omega_3$ in $\mathfrak{g}^*$.
Then 
$$d\omega_1=0,d\omega_2=\omega_1 \wedge \omega_2,
d\omega_3=-\omega_1 \wedge \omega_3.$$ So $\alpha_2=\omega_1$ and
$\alpha_3=-\omega_1$.
It is easy to see that
$H^*_{\lambda \omega}(\mathfrak{g}_0)$ is non-trivial
if and only if $\lambda\omega=0, \pm \omega_1$.
For corresponding Betti numbers
$b^p_{\lambda \omega}(M){=}
\dim H^p_{\lambda \omega}(G_0/\Gamma){=}
\dim H^p_{\lambda \omega}(\mathfrak{g_0})$ of the solvmanifold $M=G_0/\Gamma$
we have:

$$b_{\pm \omega_1}^0(M)=0, b^1_{\pm \omega_1}(M)=
b^2_{\pm \omega_1}(M)=1, b^3_{\pm \omega_1}(M)=0;$$

$$b^0(M) = b^1(M)=
b^2(M) = b^3(M) = 1.$$

\end{document}